\documentclass[11pt]{article}
\everymath{\displaystyle}
\usepackage{graphicx}
\usepackage{pst-all}
\usepackage{color}
\usepackage{amssymb}
\usepackage{amsmath}
\newtheorem{prethm}{{\bf Theorem}}

\newenvironment{thm}{\begin{prethm}{\hspace{-0.5
               em}{\bf .}}}{\end{prethm}}
\newtheorem{prelemma}{{\bf Lemma}}

\newtheorem{preex}{{\bf Example}}

\newtheorem{preprop}{{\bf Proposition}}

\newenvironment{prop}{\begin{preprop}{\hspace{-0.5em}{\bf .}}}{\end{preprop}}
\newtheorem{precor}{{\bf Corollary}}

\newenvironment{cor}{\begin{precor}{\hspace{-0.5
               em}{\bf .}}}{\end{precor}}
\newtheorem{preremark}{{\bf Remark}}

\newenvironment{remark}{\begin{preremark}{\hspace{-0.5
               em}{\bf.}}}{\end{preremark}}
\newtheorem{preprob}{{\bf Problem}}

\newenvironment{prob}{\begin{preprob}{\hspace{-0.7
               em}{\bf.}}}{\end{preprob}}
\newtheorem{predefin}{{\bf Definition}}

\newtheorem{preconj}{{\bf Conjecture}}

\newenvironment{conj}{\begin{preconj}{\hspace{-0.5
               em}{\bf .}}}{\end{preconj}}
\newtheorem{preprobb}{{\bf Problem}}

\newtheorem{prelem}{{\bf Theorem}}

\newenvironment{proof}{{\bf Proof.}\rm }{\hfill{$\Box$}}

\newtheorem{presolution}{{\bf Solution.}}

\def\newpic#1{}
\def\qed{\ifhmode\unskip\nobreak\fi\quad\ifmmode\Box\else$\Box$\fi}

\title{\vspace{-1cm}\Large\bf\noindent $(\delta, \chi_{_{\sf FF}})$-bounded families of graphs}

\author{\large\bf Manouchehr Zaker\footnote{mzaker@iasbs.ac.ir}
\vspace{5mm}\\
    Department of Mathematics,\\
     Institute for Advanced Studies in Basic Sciences,\\
    Zanjan 45137-66731, Iran\\
  }
    \date{}

\begin{document}
\maketitle

\begin{abstract}
\noindent For any graph $G$, the First-Fit (or Grundy) chromatic number of $G$, denoted by $\chi_{_{\sf FF}}(G)$, is defined as the maximum number of colors used by the First-Fit (greedy) coloring of the vertices of $G$. We call a family $\mathcal{F}$ of graphs $(\delta, \chi_{_{\sf FF}})$-bounded if there exists a function $f(x)$ with $f(x)\rightarrow \infty$ as $x\rightarrow \infty$ such that for any graph $G$ from the family one has $\chi_{_{\sf FF}}(G)\geq f(\delta(G))$, where $\delta(G)$ is the minimum degree of $G$. We first give some results concerning $(\delta, \chi_{_{\sf FF}})$-bounded families and obtain a few such families. Then we prove that for any positive integer $\ell$, $Forb(K_{\ell,\ell})$ is $(\delta, \chi_{_{\sf FF}})$-bounded, where $K_{\ell,\ell}$ is complete bipartite graph. We conjecture that if $G$ is any $C_4$-free graph then $\chi_{_{\sf FF}}(G)\geq \delta(G)+1$. We prove the validity of this conjecture for chordal graphs, complement of bipartite graphs and graphs with low minimum degree.
\end{abstract}

\noindent {\bf Mathematics Subject Classification (2000)}: 05C15, 05C07, 05C85, 05C38

\noindent {\bf Keywords:} graph coloring; First-Fit coloring; Grundy number; lower bound; minimum degree

\section{Introduction}

\noindent All graphs in this paper are simple undirected graphs. A family ${\mathcal{F}}$ of graphs is said to be $(\delta,
\chi)$-bounded if there exists a function $f(x)$ satisfying
$f(x)\rightarrow \infty$ as $x\rightarrow \infty$, such that for any
graph $G$ from the family one has $f(\delta(G))\leq \chi(G)$, where
$\delta(G)$ and $\chi(G)$ denotes the minimum degree and chromatic
number of $G$, respectively. Equivalently, the family
${\mathcal{F}}$ is $(\delta, \chi)$-bounded if $\delta(G_n)
\rightarrow \infty$ implies $\chi(G_n)\rightarrow \infty$ for any
sequence $G_1, G_2, \ldots$ with $G_n\in {\mathcal{F}}$. Motivated
by Problem 4.3 in \cite{JT}, the author introduced and
studied $(\delta, \chi)$-bounded families of graphs (under the name
of $\delta$-bounded families) in \cite{Z}. The so-called color-bound
family of graphs mentioned in the related problem of \cite{JT} is a
family for which there exists a function $f(x)$ satisfying
$f(x)\rightarrow \infty$ as $x\rightarrow \infty$, such that for any
graph $G$ from the family one has $f(col(G))\leq \chi(G)$, where
$col(G)$ is defined as $col(G)=\max \{\delta(H): H\subseteq G\}+1$. It was shown in \cite{Z} that if we restrict ourselves to hereditary (i.e.
closed under taking induced subgraph) families then the two concepts
$(\delta, \chi)$-bounded and color-bound are equivalent. The first
specific results concerning $(\delta, \chi)$-bounded families
appeared in \cite{Z} where the following theorem was proved (in a
somewhat different but equivalent form). In the following theorem for any set ${\mathcal{C}}$
of graphs, $Forb({\mathcal{C}})$ denotes the class of graphs that
contains no member of ${\mathcal{C}}$ as an induced subgraph.

\begin{thm} (\cite{Z})
For any set ${\mathcal{C}}$ of graphs, $Forb({\mathcal{C}})$ is
$(\delta, \chi)$-bounded if and only if there exists a constant
$c=c({\mathcal{C}})$ such that for any bipartite graph $H\in
Forb({\mathcal{C}})$ one has $\delta(H)\leq c$.
\label{turan}
\end{thm}

\noindent Theorem \ref{turan} shows that to decide whether $Forb({\mathcal{C}})$ is $(\delta,
\chi)$-bounded we may restrict ourselves to bipartite graphs. A comprehensive study of $(\delta,
\chi)$-bounded families was done in \cite{GZ}, where the authors proved the following
theorem.

\begin{thm}(\cite{GZ})
Given a finite set of graphs $\{H_1, H_2, \ldots, H_k\}$. Then $Forb(H_1,
H_2, \ldots, H_k)$ is $(\delta, \chi)$-bounded if and only if one of the
following holds:

\noindent (i) For some $i$, $H_i$ is a star tree.

\noindent (ii) For some $i$, $H_i$ is a forest and for some $j\not= i$, $H_j$ is
complete bipartite graph.\label{finitecollect}
\end{thm}

\noindent The following result concerns $Forb({\mathcal{C}})$, where ${\mathcal{C}}$ contains infinitely many graphs in which one of them is a tree.

\begin{thm}(\cite{GZ})
Let $T$ be any non star tree. Then $Forb(T, H_1, \ldots)$ is
$(\delta, \chi)$-bounded if and only if at least one of $H_i$-s is
complete bipartite graph.
\end{thm}

\noindent For other $(\delta, \chi)$-bounded families of graphs we refer the reader to \cite{GZ} and \cite{Z}. In this paper we work on the Grundy (or First-Fit) chromatic number of graphs. A {\it Grundy $k$-coloring} of a graph $G$, is a proper
$k$-coloring of vertices in $G$ using colors $\{1, 2, \ldots, k\}$
such that for any two colors $i$ and $j$ with $i<j$, any vertex
colored $j$ is adjacent to some vertex colored $i$. The Grundy or
{\it First-Fit chromatic number} of a graph $G$, denoted by
$\chi_{{_{\sf FF}}}(G)$ (also denoted by $\Gamma(G)$ in some articles),
is the largest integer $k$, such that there exists a Grundy
$k$-coloring for $G$. It can be shown that $\chi_{{_{\sf FF}}}(G)$ is the same as the maximum number of colors used by the First-Fit (greedy) coloring of the vertices of $G$ \cite{Z2}. To determine $\chi_{{_{\sf FF}}}(G)$ is NP-complete even for complement of bipartite graphs $G$ \cite{Z2}. For this reason it is natural to obtain lower and upper bounds for $\chi_{{_{\sf FF}}}(G)$ in terms of ordinary graph theoretical parameters. In this paper we obtain some lower bounds in terms of the minimum degree of graphs. The Grundy number and First-Fit coloring
of graphs were studied widely in the literature, see \cite{Z2,Z1} and
its references. Throughout the paper we denote the complete graph on $n$ vertices by $K_n$ and the cycle on $n$ vertices by $C_n$. For each positive integer $\ell$, the complete bipartite graph in which each part has $\ell$ vertices is denoted by $K_{\ell,\ell}$.

\section{Results}

\noindent Generalizing the concept of $(\delta, \chi)$-bounded graph, we define the following notion. A family $\mathcal{F}$ of graphs is called $(\delta, \chi_{_{\sf FF}})$-bounded if there exists a function $f(x)$ with $f(x)\rightarrow \infty$ as $x\rightarrow \infty$ such that for any graph $G$ from the family one has $\chi_{_{\sf FF}}(G)\geq f(\delta(G))$. It was shown in \cite{Z2} that $\chi_{_{\sf FF}}(G)=2$ if and only if $G$ is a complete bipartite graph. Obviously, complete bipartite graphs may have arbitrary large minimum degree. We conclude that the family of complete bipartite graphs is not $(\delta, \chi_{_{\sf FF}})$-bounded. This example induces that perhaps $C_4$ and other complete bipartite graphs have significant role in $(\delta, \chi_{_{\sf FF}})$-boundedness. Note also that any $(\delta, \chi)$-bounded family is also $(\delta, \chi_{_{\sf FF}})$-bounded. Another interesting chromatic-related parameter is the so-called coloring number of graphs. The coloring number of a graph $G$ is defined as $col(G)={\max}_{H\subseteq G} \delta(H)+1$. The coloring number of graphs is a polynomial time parameter. See \cite{JT,MGR, Z} for more results on the coloring number of graphs. The possible relationships between the coloring number and Grundy number of graphs is an interesting research area. For some graphs $G$ we have $col(G) < \chi_{_{\sf FF}}(G)$. For example the path on four vertices $P_4$ and infinitely many trees satisfy this inequality. Also, for some graphs $G$ we have $\chi_{_{\sf FF}}(G) < col(G)$. For example consider complete bipartite graphs $K_{a,b}$, where $a,b\geq 2$. We have the following remark, where by a hereditary family we mean any family ${\mathcal{F}}$ such that for any graph $G$ from ${\mathcal{F}}$, if $H$ is an induced subgraph of $G$ then $H\in {\mathcal{F}}$.

\begin{remark}
Let ${\mathcal{F}}$ be any hereditary family of graphs such that ${\mathcal{F}}$ is $(\delta, \chi_{_{\sf FF}})$-bounded. Then there exists a function $f(x)$ with $f(x)\rightarrow \infty$ as $x\rightarrow \infty$ such that for any graph $G$ from the family one has $\chi_{_{\sf FF}}(G)\geq f(col(G))$.
\end{remark}

\noindent \begin{proof}
Let $G \in {\mathcal{F}}$ and $H$ be any induced subgraph of $G$ with $col(G)=\delta(H)+1$. We have $H\in {\mathcal{F}}$. Let $g(x)$ be such that
$\chi_{_{\sf FF}}(H)\geq g(\delta(H))$. The proof completes by taking $f(x)=g(x-1)$.
\end{proof}

\noindent As we mentioned before any $(\delta, \chi)$-bounded family is also $(\delta, \chi_{_{\sf FF}})$-bounded. In Theorem \ref{log} we obtain $(\delta, \chi_{_{\sf FF}})$-bounded families which are not $(\delta, \chi)$-bounded. For this purpose we first obtain in Proposition \ref{girth} a result concerning $(\delta, \chi)$-boundedness of graphs. In the following, the girth of a graph $G$ is the length of a shortest cycle contained in $G$. When $G$ contains a cycle then we say that $G$ has finite girth. We use also the following two facts. The first fact states that any graph with $m$ edges contains a bipartite subgraph with at least $m/2$
edges. The second one states that any graph with $n$ vertices and $m/2$ edges contains a subgraph
with minimum degree at least $m/2n$. For the proof of these facts we refer the reader to standard Graph Theory books such as \cite{D}.

\begin{prop}
Let ${\mathcal{C}}$ be any finite collection of graphs such that any member of ${\mathcal{C}}$ has finite girth. Then $Forb({\mathcal{C}})$ is not $(\delta, \chi)$-bounded. In particular $Forb(K_3, K_{2,m})$ and $Forb(K_{\ell,\ell})$ are not $(\delta, \chi)$-bounded.\label{girth}
\end{prop}

\noindent \begin{proof}
Let $g$ be an even integer such that the girth of any graph in ${\mathcal{C}}$ is at most $g$.
For the proof we use
the following Tur\'an-type result which is attributed to Erd\H{o}s in
\cite{S}. For any $k$ and $n$ there exists a graph on $n$ vertices
with $\Omega(n^{1+1/2k-1})$ edges that contains no cycle of length at most $2k$. Let $g=2k$, and recall from the previous paragraph that (1)
a graph with $m$ edges contains a bipartite subgraph with at least $m/2$
edges, and (2) a graph with $n$ vertices and $m/2$ edges contains a subgraph
with minimum degree at least $m/2n$. We conclude that there exists an infinite sequence $G_1, G_2, \ldots$ of bipartite
graphs such that $\delta(G_i)\rightarrow \infty$ as $i\rightarrow
\infty$ and the girth of any $G_i$ is more than $g$. This shows that
$G_i$ belongs to $Forb({\mathcal{C}})$. This shows that $Forb({\mathcal{C}})$ is not $(\delta, \chi)$-bounded.
\end{proof}

\noindent In opposite direction we show in Theorem \ref{log} that $Forb(K_3, K_{2,m})$ is
$(\delta, \chi_{_{\sf FF}})$-bounded. More generally, Theorem \ref{referee} asserts that $Forb(K_{\ell,\ell})$ is
$(\delta, \chi_{_{\sf FF}})$-bounded. Before we proceed, we need to introduce a family of trees $T_k$, $k=1, 2, \ldots$. For $k=1,2$, $T_1$ (resp. $T_2$) is isomorphic to $K_1$ (resp. $K_2$). Assume that $T_k$ has been defined. Attach a leaf to any vertex of $T_k$ and denote the resulting tree by $T_{k+1}$. It is easily observed that $\chi_{_{\sf FF}}(T_k)=k$. Note also that $|V(T_k)|=2^{k-1}$. We need also the following result from (\cite{Z}, Theorem 2), where $\rho(G)=|E(G)|/|V(G)|$.

\begin{thm}(\cite{Z})
Let $G$ be any triangle-free graph such that $G$ does not contain $K_{2,m}$, where $m>1$. If $\rho(G)\geq (k-3)(m-1)+1$ then $G$ contains
all trees on $k$ vertices as induced subgraphs.\label{rho}
\end{thm}

\noindent The promised result is as follows.

\begin{thm}
Let $G\in Forb(K_3, K_{2,m})$. Then $$\chi_{_{\sf FF}}(G)\geq \log (\frac{\delta(G)+6m-8}{2m-2})+1.$$\label{log}
\end{thm}

\noindent \begin{proof}
First note that $G$ does not contain triangle and $K_{2,m}$. Set $\delta(G)=p$ and $k=(p+6m-8)/(2m-2)$, for simplicity.
We have
$$(k-3)(m-1)+1=(\frac{p+6m-8}{2m-2}-3)(m-1)+1$$
$$\hspace{2cm}=(\frac{p-2}{2m-2})(m-1)+1$$
$$\hspace{-1.1cm}=\frac{p}{2}.$$
\noindent Since $\rho(G)\geq (p/2)$ then $G$ satisfies the conditions of Theorem \ref{log} with these $k$ and $m$. Therefore $G$ contains all trees on $k$ vertices as induced subgraph. In particular $G$ contains $T_q$ as induced subgraph, where $q= \log ((\delta(G)+6m-8)/(2m-2))+1$. We conclude that $$\chi_{_{\sf FF}}(G)\geq \log (\frac{\delta(G)+6m-8}{2m-2})+1.$$
\end{proof}

\noindent By applying Theorem \ref{rho} and Theorem \ref{log} more economically when $m=2$ we obtain the following bound.

\begin{cor}
Let $G\in Forb(K_3, C_4)$. Then $$\chi_{_{\sf FF}}(G)\geq \log (\delta(G)+1).$$\label{3,4}
\end{cor}

\noindent We noted before that the family of complete bipartite graphs is not $(\delta, \chi_{_{\sf FF}})$-bounded. Hence the following proposition is immediate from this fact.

\begin{prop}
Let ${\mathcal{C}}$ be any collection of graphs such that any member of it contains an odd cycle. Then $Forb({\mathcal{C}})$ is not $(\delta, \chi_{_{\sf FF}})$-bounded. In particular $Forb(K_3)$ is not $(\delta, \chi_{_{\sf FF}})$-bounded.
\end{prop}

\noindent In Theorem \ref{referee} we prove that $Forb(K_{\ell,\ell})$ is $(\delta, \chi_{_{\sf FF}})$-bounded. For this purpose we need the following theorem from \cite{GZ}.

\begin{thm}(\cite{GZ})
For every tree $T$ and for positive integers $\ell$, $k$ there exist a function $f(T, \ell, k)$
with the following property. If $G$ is a graph with $\delta(G) \geq f(T, \ell, k)$ and $\chi(G) \leq k$ then $G$
contains either $T$ or $K_{\ell,\ell}$ as an induced subgraph.\label{GZ}
\end{thm}

\noindent We shall make use of this theorem in proving the next result.

\begin{thm}
For each positive integer $\ell$,  $Forb(K_{\ell,\ell})$ is $(\delta, \chi_{_{\sf FF}})$-bounded.\label{referee}
\end{thm}

\noindent \begin{proof}
Recall that for each positive integer $k$, $T_k$ denotes the only smallest tree of Grundy number $k$. Let $\{G_n\}_{n=1}^{\infty}$ be a sequence of $K_{\ell,\ell}$-free graphs such that $\delta(G_n)\rightarrow \infty$ as $n\rightarrow \infty$. Assume on the contrary that for some integer $N$, $\chi_{_{\sf FF}}(G_n)\leq N$ holds for all $n$. It follows that for each $n$, $T_{N+1}$ is not an induced subgraph of $G_n$. Hence $G_n$ belongs to $Forb(T_{N+1},K_{\ell, \ell})$. Theorem \ref{GZ} implies that for each $n$ either $\delta(G_n)< f(T_{N+1}, \ell, N)$ or $\chi(G_n)>N$. But the second case is impossible because $\chi(G_n)\leq \chi_{_{\sf FF}}(G_n)\leq N$. Therefore for each $n$, $\delta(G_n)< f(T_{N+1}, \ell, N)$. This contradicts $\delta(G_n)\rightarrow \infty$. This contradiction completes the proof.
\end{proof}

\noindent The following result shows that chordal graphs are $(\delta, \chi_{_{\sf FF}})$-bounded with $f(x)=x+1$. Note that the class of chordal graphs is the same as $Forb(C_4, C_5, \ldots)$. In a graph $G$, by a simplicial vertex we mean any vertex $v$ such that $G[N(v)]$ is a clique in $G$, where $G[N(v)]$ stands for the subgraph of $G$ induced by the set $N(v)$ of the neighbors of $v$ in $G$. It is a known fact that any chordal graph $G$ admits a simplicial elimination ordering (see e.g. \cite{D}). In other words, let $G$ be a chordal graph. Then there exists a vertex ordering $v_1, \ldots, v_n$ of $G$ such that $v_i$ is simplicial in $G\setminus \{v_1, \ldots, v_{i-1}\}$. We shall make use of this fact in the following theorem.

\begin{thm}
$Forb(C_4, C_5, \ldots)$ is $(\delta, \chi_{_{\sf FF}})$-bounded with $f(x)=x+1$.\label{x+1}
\end{thm}

\noindent \begin{proof}
Let $G$ be any chordal graph $G$ and let $v_1, v_2, \ldots, v_n$ be a simplicial ordering of the vertices of $G$. Since $G[N(v)\cup \{v_1\}]$ is a clique in $G$ with $deg_G(v_1)+1$ vertices, then $\omega(G) \geq deg_G(v) +1 \geq \delta(G)+1$. We have also $\chi_{_{\sf FF}}(G)\geq \chi(G) \geq \omega(G)$. Hence $\chi_{_{\sf FF}}(G) \geq \delta(G)+1$.
\end{proof}

\noindent By strengthening Theorem \ref{x+1} we propose the following conjecture.

\begin{conj}
Let $G$ be a $C_4$-free graph. Then $\chi_{_{\sf FF}}(G)\geq \delta(G)+1$.\label{conj}
\end{conj}

\noindent A natural scenario to prove the above conjecture is as follows. Let ${\mathcal{F}}$ be a hereditary family of graphs satisfying the
following property. Any member $G$ from the family contains a maximal independent set (MIS) such as $I$ such that $\delta(G\setminus I)=\delta(G)-1$.
We have the following observation which can be proved by induction.\\

\noindent {\bf Observation 1.}
{\it Let ${\mathcal{F}}$ be any hereditary family of graphs such that any graph $G$ from the family contains a MIS, say $I$ such that for any vertex $v$ of $G$ if $deg_G(v)=\delta(G)$ then $deg_{G\setminus I}(v)=deg_G(v)-1$. Then $\chi_{_{\sf FF}}(G)\geq \delta(G)+1$ for any graph $G$ from ${\mathcal{F}}$.}\\

\noindent Unfortunately
the family of $C_4$-free graphs does not satisfy the above condition. In this regard it is worthy to work on the following problem.

\begin{prob}
Find families ${\mathcal{F}}$ of graphs satisfying the following property. Any graph $G$ from ${\mathcal{F}}$ contains a MIS, say $I$ such that for any vertex $v$ of $G$ if $deg_G(v)=\delta(G)$ then $deg_{G\setminus I}(v)=deg_G(v)-1$.
\end{prob}

\noindent In Theorem \ref{cobip} we show that Conjecture \ref{conj} holds for any graph which is the complement of a bipartite graph. We need some prerequisites. In a graph $H$, a subset $D$ of edges in $H$ is called an edge dominating set if each edge
in $E(H)\setminus D$ has a common end point with an edge in $D$. Let $H$ be any bipartite graph. Set $G=\overline{H}$. Let $\gamma'(H)$ be the smallest size of an edge dominating set in $H$. It was proved in \cite{Z2} that $\chi_{_{\sf FF}}(G)=|V(G)|-\gamma'(H)$. We have now the following theorem.

\begin{thm}
Let $H$ be any bipartite graph and $G$ be the complement of $H$ such that $G$ is $C_4$-free. Then $\chi_{_{\sf FF}}(G)\geq \delta(G)+1$.\label{cobip}
\end{thm}

\noindent \begin{proof}
Let $n$ be the order of $G$. Since $\delta(G)=n-\Delta(H)-1$, then the inequality $\chi_{_{\sf FF}}(G)\geq \delta(G)+1$ is equivalent to $\gamma'(H)\leq \Delta(H)$. Now we use the fact that for any edge dominating set $R$ in a bipartite graph, there is a matching $M$ which is also an edge dominating set and $|M|\leq |R|$. This fact can be easily proved and we omit mentioning its proof here, and refer the reader to \cite{HAD}. Let $R$ be an edge dominating set in $H$ with $|R|=\gamma'(H)$. Using the previous fact we obtain that $R$ is a matching and therefore $\gamma'(H)\leq \alpha'(H)$. Hence to complete the proof we need to show that $\alpha'(H)\leq \Delta(H)$. We prove the latter inequality by induction on the number of edges. Note that since $G$ is $C_4$-free then $H$ is $2K_2$-free, where by $2K_2$ we mean the graph consisting only of two independent edges.

\noindent Let $M$ be any matching of maximum size in $H$ and $e=uv$ be any edge of $M$. Define another graph as $H_0=H\setminus \{u, v\}$. We have $\alpha'(H_0)\leq \Delta(H_0)$. Hence $\alpha'(H)-1\leq \Delta(H_0)$. To finalize the proof we show that $\Delta(H_0)+1\leq \Delta(H)$. Otherwise, since $H_0$ is an induced subgraph of $H$, we have $\Delta(H_0)= \Delta(H)$. Let $x$ be any vertex in $H_0$ such that $deg_{H_0}(x)=\Delta(H)$. Without loss of generality, we may assume that $u$ and $x$ are in the same bipartite part of $H$. We show that $u$ is adjacent to any neighbor of $x$. Let $w$ be any neighbor of $x$. Since $H$ is $2K_2$-free then the subgraph of $H$ consisting of two edges $uv$ and $xw$ can not be induced in $H$. Now, since $x$ has the maximum degree then $x$ can not be adjacent to $v$ in $H$. Hence $u$ should be adjacent to $w$. But $v$ is adjacent to $u$ and not adjacent to $x$. This means that the degree of $u$ is strictly greater that the degree of $x$, a contradiction with our choice of $x$. This completes the proof.
\end{proof}

\noindent The following theorem shows that Conjecture \ref{conj} holds for all graphs $G$ with $\delta(G)\leq 3$.

\begin{thm}
Let $G$ be a $C_4$-free graph with $\delta(G)\leq 3$. Then $\chi_{_{\sf FF}}(G)\geq \delta(G)+1$.
\end{thm}

\noindent \begin{proof}
Theorem obviously holds if $\delta(G)=1$. Assume that $\delta(G)=2$. Let $v$ be any vertex of $G$ and $a, b$ any two neighbors of $v$. If $a$ and $b$ are adjacent then the resulting triangle shows that $\chi_{_{\sf FF}}(G)\geq 3$. Assume that $a$ and $b$ are not adjacent. Let $c$ be any neighbor of $b$. If $a$ and $c$ are not adjacent then we obtain an induced $P_4$ on the vertex set $\{v, a, b, c\}$. Hence the desired inequality holds in this case. Assume that $a$ and $c$ are adjacent. Then since $G$ is $C_4$-free and $b$ is not adjacent to $a$ then $v$ is adjacent to $c$. This gives rise to a triangle. Hence in this case too $\chi_{_{\sf FF}}(G)\geq 3$.

\noindent Assume now that $\delta(G)=3$. Let $v$ be any vertex of degree 3. Let $a, b, c$ be the neighbors of $v$.

\noindent {\bf Case 1.} Three vertices $a, b, c$ are independent.

\noindent In this case we first note that no two vertices from $\{a, b, c\}$ have a common neighbor other than the vertex $v$, since otherwise let $u$ be a common neighbor of $a$ and $b$. The two vertices $a$ and $b$ are independent and $G$ is $C_4$-free. Hence $v$ should be adjacent to $u$. This contradicts $deg_G(v)=3$. Now, let $x$ and $y$ (resp. $z$ and $t$) be two neighbors of $a$ (resp. $b$). We have $\{x,y\}\cap \{z,t\}=\varnothing$ and $c$ is not adjacent to any vertex in $\{x, y, z, t\}$. Consider a small bipartite graph $H$ consisting of the bipartite sets $\{x, y\}$ and $\{z, t\}$ with all edges from $G$ among these parts. If there are at most two edges in $H$ then we color $v$ by 4, $a$ by 3, $b$ by 2 and $c$ by 1. The vertex $a$ needs two neighbors of colors 1 and 2; and the vertex $b$ needs one neighbor of color 1. We can easily fulfil these conditions by assigning suitable colors 1 and 2 to the vertices of $H$. If there are exactly three edges in $H$ then (assuming that $x$ is adjacent to $z$) we consider the following coloring. We color $x$ by 4, $z$ by 3, $y$ and $t$ by 2 and $a$ and $b$ by 1. Finally, we consider the case that $H$ is a complete bipartite graph. In this case we color $x$ by 4, $z$ by 3, $a$ and $b$ by 2; and $t$ and $v$ by 1 (note that $t$ and $v$ are not adjacent). All of these pre-colorings are partial Grundy colorings with 4 colors. This completes the proof in Case 1.

\noindent {\bf Case 2.} $a$ and $b$ are adjacent and $c$ is not adjacent to $a$ and also to $b$.

\noindent In this case we color $v$ by 4, $a$ by 2, $b$ by 3 and $c$ by 1.
Let $d$ be a neighbor of $b$. We color $d$ by 1. Note that $c$ and $d$ can not be adjacent. If $a$ is adjacent to $d$ then we obtain a partial Grundy coloring using four colors. Otherwise, let $a$ be adjacent to a vertex say $e$. If $e$ is adjacent to $d$ then since $a$ is not adjacent to $d$, hence $b$ should be adjacent to $e$. In this situation we color $e$ by 1 and remove the color of $d$. Now the colors of $\{v, a, b, c, e\}$ is a partial Grundy coloring with four colors. But if $e$ is not adjacent to $d$, we color both vertices $e$ and $d$ by 1. Note that in this case the colors of $\{v, a, b, c, d, e\}$ introduce a partial Grundy coloring using four colors.

\noindent {\bf Case 3.} $a$ is adjacent to both $b$ and $c$; but $b$ and $c$ are not adjacent.

\noindent In this case we color $v$ by 4, $a$ by 3, $b$ by 1 and $c$ by 2.
Let $d$ be a new neighbor of $c$. If $b$ and $d$ are not adjacent then we color both of them by 1. The resulting coloring is a partial Grundy coloring using four colors. But if $b$ and $d$ are adjacent, then $a$ should be adjacent to $d$. Now consider the 4-cycle on $\{v, b, c, d\}$. Since the degree of $v$ is three then $v$ can not be adjacent to $d$. Hence $b$ and $c$ should be adjacent. But this is a contradiction.

\noindent {\bf Case 4.} The only remaining case is that $a, b, c$ are all adjacent. But in this case we obtain a clique of size four. It is clear that in this case $\chi_{_{\sf FF}}(G)\geq 4$.
\end{proof}

\noindent We end the paper by mentioning that Conjecture \ref{conj} is also valid for graphs with minimum degree four. The proof is by checking too many cases. We omit the details.

\end{document}